\documentclass{amsart}
\usepackage{amssymb}
\usepackage{amscd}

\newcommand\Rs{{\mathbb R}}
\newcommand\ks{{\Bbbk}}

\newcommand\Der{\mathop {\fam 0 Der}\nolimits}
\newcommand\id{\mathop {\fam 0 id}\nolimits}

\begin{document}

\title{Commutation relations on the covariant derivative}
\thanks{Supported in part by the grant RFFI-08-01-92001
}
\author{Alexey~V.~Gavrilov}
\email{gavrilov19@gmail.com}

\maketitle

Let ${\mathfrak g}$ be a nonassociative algebra over a field $\ks.$
The operation in ${\mathfrak g}$ will be denoted by the diamond sign $\diamond$,
for example,
$(x\diamond y)\diamond z\in {\mathfrak g},\,x,y,z\in{\mathfrak g}.$
Let $T({\mathfrak g})=T_{\ks}({\mathfrak g})=
\oplus_{n=0}^{\infty}{\mathfrak g}^{\otimes n}$ denotes the tensor algebra of
${\mathfrak g}$ and
$L_u:T({\mathfrak g})\to T({\mathfrak g}),\,u\in T({\mathfrak g})$, denotes
the operator of left multiplication, $L_u:v\mapsto u\otimes v.$
For typographical reasons we shall write sometimes $L(u)$ instead of $L_u.$
Denote by $\tau_x\in\Der_{\ks}(T({\mathfrak g})),\,x\in{\mathfrak g}$
the derivation of the tensor algebra defined by the condition
$\tau_x:y\mapsto x\diamond y,\,y\in{\mathfrak g}.$
There exists a unique linear map
$K:T({\mathfrak g})\to T({\mathfrak g})$ such that
$K(1)=1$ and $KL_x=L_xK-K\tau_x,\,x\in{\mathfrak g}.$
For example, $K(x)=K(L_x 1)=L_x(K(1))-K(\tau_x 1)=x,
\,K(x\otimes y)=K(L_x y)=L_x(K(y))-K(\tau_x y)=
x\otimes y-x\diamond y,$ etc.
The author introduced this map in [1,2] (it is
unlikely that such a map has never been considered before, however,
I have not found an appropriate reference).

The problem we are interested in appears when ${\mathfrak g}$ is
a nonassociative algebra and a Lie algebra simultaneously
(a "framed Lie algebra" in terms of [1]). In other words,
it is an algebra with two operations, one of which is antisymmetric and
satisfies the Jacobi identity.
As ${\mathfrak g}$ is a Lie algebra, there exists the exact sequence
$0\to I({\mathfrak g})\to T({\mathfrak g})\to U({\mathfrak g})\to 0,$
where $U({\mathfrak g})$ is the universal enveloping algebra of ${\mathfrak g}$
and  $I({\mathfrak g})$ is the two-sided ideal of $T({\mathfrak g})$
generated by the elements of the form $x\otimes y-y\otimes x-[x,y].$
Let $\Sigma({\mathfrak g})\subset T({\mathfrak g})$ denotes the linear
space of symmetric tensors (i.e. $\Sigma_n({\mathfrak g})=
\Sigma({\mathfrak g})\cap{\mathfrak g}^{\otimes n}=
({\mathfrak g}^{\otimes n})^{S_n}$).
The map $K$ has the property $\deg(Ku-u)<\deg(u),\,u\in T({\mathfrak g})$
hence it is invertible. Moreover, if ${\rm char}\ks=0,$
then the restriction of the natural map to
$K(\Sigma({\mathfrak g}))\to U({\mathfrak g})$
is a linear isomorphism by the Poincare-Birkhoff-Witt theorem.
Thus we have the decomposition
$T({\mathfrak g})=\Sigma({\mathfrak g})\oplus K^{-1}(I({\mathfrak g})).$
The problem is to describe the corresponding projection
$\pi:T({\mathfrak g})\to\Sigma({\mathfrak g})$  explicitly.

The origin of the problem lies in differential geometry; it will be discussed
below. Even in the case $K=\id$ (i.e.
${\mathfrak g}\diamond{\mathfrak g}=\{0\}$) it is not trivial
[5]. In the general case it looks like a formidable task.
This paper concerned with a much more simple related problem. Denote
$\Omega({\mathfrak g})=T({\mathfrak g})\otimes\wedge^2{\mathfrak g},$
where $\wedge^2{\mathfrak g}$ is the exterior square of ${\mathfrak g}$.
Let $\imath:\Omega({\mathfrak g})\otimes T({\mathfrak g})
\to T({\mathfrak g})$
be the natural map, $\imath:a\otimes x\wedge y\otimes b
\mapsto a\otimes(x\otimes y-y\otimes x)\otimes b$.
The aim is to find in an explicit form a map
$R:\Omega({\mathfrak g})\otimes T({\mathfrak g})\to T({\mathfrak g})$
satisfying for any $Q\in \Omega({\mathfrak g})\otimes T({\mathfrak g})$
the following two properties: $R(Q)+\imath(Q)\in K^{-1}(I({\mathfrak g}))$
and $\deg R(Q)<\deg\imath(Q).$
These properties do not determine the map $R$ uniquely;
however, the solution proposed below is probably the simplest one.

The projection $\pi$ may be expressed via $R$ as follows. Let
$s:T({\mathfrak g})\to S({\mathfrak g})$ be the natural algebra homomorphism
and $J({\mathfrak g})=\ker(s,T({\mathfrak g})).$ The restriction
$s:\Sigma({\mathfrak g})\to S({\mathfrak g})$ is a linear isomorphism, hence
$T({\mathfrak g})=\Sigma({\mathfrak g})\oplus J({\mathfrak g}).$
Denote by $\pi_0:T({\mathfrak g})\to \Sigma({\mathfrak g})$ the corresponding
projection. The map
$\imath:\Omega({\mathfrak g})\otimes T({\mathfrak g})\to J({\mathfrak g})$ is
surjective, hence there exists a (non-unique) right inverse
$\imath^{-1}:J({\mathfrak g})\to \Omega({\mathfrak g})\otimes T({\mathfrak g}).$
Denote $\rho=
R\circ\imath^{-1}\circ (\id-\pi_0):T({\mathfrak g})\to T({\mathfrak g}).$
One can choose an inverse map $\imath^{-1}$ satisfying
the natural condition $\deg\imath^{-1}(u)<\deg u,\,u\in T({\mathfrak g})$.
Then we have $\deg\rho(u)<\deg u$, hence
the map $\id+\rho$ is invertible. Let
$\pi=\pi_0(\id+\rho)^{-1}.$ By the definitions,
$\rho:\Sigma({\mathfrak g})\to\{0\}$ and
$\id+\rho:J({\mathfrak g})\to K^{-1}(I({\mathfrak g})).$
Hence $\pi$ is exactly the projection we need (though not in an explicit form).
We shall see below that the map $R$ may be
interpreted as a collection of commutation
relations on the covariant derivative.

\section{ Relations}

In this section, ${\mathfrak g}$ is an algebra with two operations:
the first one denoted by $\diamond$ and the antisymmetric
second one denoted by the brackets
$[\cdot,\cdot].$ In the applications it is a Lie algebra,
i.e. the second operation satisfies the Jacobi identity,
but we shall not actually need this identity here.
The base field $\ks$ is an arbitrary one.

Let $I({\mathfrak g})\subset T({\mathfrak g})$ be the two-sided
ideal generated by the
elements of the form $x\otimes y-y\otimes x-[x,y],\,x,y\in{\mathfrak g}.$
Note that if ${\mathfrak g}$ is not a Lie algebra then
$T({\mathfrak g})/I({\mathfrak g})$ is no more the
universal enveloping algebra.
Let $\Omega({\mathfrak g})$ and $\tau_x,\,x\in{\mathfrak g}$
are defined as above. The derivation
$\tau_x:T({\mathfrak g})\to T({\mathfrak g})$
can be naturally lifted to the map
$\tau_x:\Omega({\mathfrak g})\to\Omega({\mathfrak g})$
by the condition $\tau_x\circ\imath=\imath\circ\tau_x;$
namely $\tau_x:a\otimes y\wedge z\mapsto\tau_x(a)\otimes y\wedge z+
a\otimes x\diamond y\wedge z+a\otimes y\wedge x\diamond z.$
The linear maps $t:\Omega({\mathfrak g})\to{\mathfrak g},
r:\Omega({\mathfrak g})\to\Der_{\Bbbk}(T({\mathfrak g}))$
and $e:\Omega({\mathfrak g})\to I({\mathfrak g})$ are defined as follows.
If $x,y,z\in{\mathfrak g}$ and $Q\in\Omega({\mathfrak g}),$ then
$$t(x\wedge y)=x\diamond y-y\diamond x-[x,y],\,
t(x\otimes Q+\tau_x Q)= x\diamond t(Q),$$
$$r(x\wedge y):z\mapsto
x\diamond(y\diamond z)-y\diamond(x\diamond z)-[x,y]\diamond z,\,
r(x\otimes Q+\tau_x Q)=[\tau_x,r(Q)],$$
$$e(x\wedge y)=x\otimes y-y\otimes x-[x,y],\,
e(x\otimes Q+\tau_x Q)= x\otimes e(Q)-e(Q)\otimes x.$$

It is well known that tensor algebra is a bialgebra
(actually it is a Hopf algebra but we make no use of antipode).
The comultiplication $\bigtriangleup:T({\mathfrak g})\to
T({\mathfrak g})\otimes T({\mathfrak g})$ is an algebra homomorphism
defined by the condition $\bigtriangleup:
x\mapsto 1\hat{\otimes} x+x\hat{\otimes} 1,\,x\in{\mathfrak g}$
(to avoid misunderstanding the "exterior" tensor product
is denoted by the hatted sign,
i.e. $1\hat{\otimes} x$ and $x\hat{\otimes} 1$ are the elements of
$T({\mathfrak g})\otimes T({\mathfrak g})$, not of $T({\mathfrak g})$).
We shall use the common Sweedler notation, e.g.
$\bigtriangleup(u)=\sum_{(u)}u_{(1)}\hat{\otimes} u_{(2)}.$

{\bf Theorem. }
{\it Let ${\mathfrak g}$ be as above,
$u,v\in T({\mathfrak g})$ and $\omega\in\wedge^2{\mathfrak g}.$
Then
$$K(u\otimes\imath(\omega)\otimes v+
\sum_{(u)}u_{(1)}\otimes
(t(u_{(2)}\otimes\omega)\otimes v+r(u_{(2)}\otimes\omega)v))=
\sum_{(u)}e(u_{(1)}\otimes\omega)\otimes K(u_{(2)}\otimes v).$$
}

As an easy consequence, the map
$$R:\Omega({\mathfrak g})\otimes T({\mathfrak g})\to T({\mathfrak g}),\,
R:u\otimes\omega\otimes v\mapsto
\sum_{(u)}u_{(1)}\otimes
(t(u_{(2)}\otimes\omega)\otimes v+r(u_{(2)}\otimes\omega)v)$$
has the required properties: if
$Q\in\Omega({\mathfrak g})\otimes T({\mathfrak g}),$ then
$R(Q)+\imath(Q)\in K^{-1}(I({\mathfrak g}))$ and
$\deg R(Q)<\deg\imath(Q).$

It is convenient to introduce the linear maps
$\lambda_x=L_x+\tau_x,\,x\in{\mathfrak g}$
and $q(Q)=L_{t(Q)}+r(Q):T({\mathfrak g})
\to T({\mathfrak g}),\,Q\in\Omega({\mathfrak g}).$

Denote by
$Z(u,\omega):T({\mathfrak g})\to T({\mathfrak g}),
\,u\in T({\mathfrak g}),\omega\in\wedge^2{\mathfrak g}$
the linear map defined by the equality
$$Z(u,\omega)=KL(u)L(\imath(\omega))+
\sum_{(u)}KL(u_{(1)})q(u_{(2)}\otimes\omega)-
L(e(u_{(1)}\otimes\omega))KL(u_{(2)}).$$
The statement of the theorem may be written in the form
$Z(u,\omega)v=0,$ so it remains to prove that $Z(u,\omega)$
is zero.
By the definitions,
$$Z(1,x\wedge y)=K[L_x,L_y]+Kq(x\wedge y)-([L_x,L_y]-L_{[x,y]})K.$$
Substituting
$$q(x\wedge y)=[\tau_x,L_y]-[\tau_y,L_x]+[\tau_x,\tau_y]-\lambda_{[x,y]}$$
and taking into account the equality $K\lambda_{x}=L_xK,\,x\in{\mathfrak g},$
we get
$$Z(1,x\wedge y)=K[\lambda_x,\lambda_y]-[L_x,L_y]K-K\lambda_{[x,y]}
+L_{[x,y]}K=0.$$

An easy computation shows that
$L(\lambda_x u)=\lambda_x L(u)-L(u)\tau_x,
q(\lambda_x Q)=[\tau_x, q(Q)]$ and
$\bigtriangleup(\lambda_xu)=\sum_{(u)}\lambda_xu_{(1)}\hat{\otimes} u_{(2)}+
u_{(1)}\hat{\otimes} \lambda_xu_{(2)}$
for any $x\in{\mathfrak g}$.
Applying all these equalities we have
$$Z(\lambda_x u,\omega)+Z(u,\tau_x\omega)=
KL(\lambda_x(u\otimes\imath(\omega)))+
KL(\lambda_x u_{(1)})q(u_{(2)}\otimes\omega)+$$
$$+KL(u_{(1)})q(\lambda_x(u_{(2)}\otimes\omega))-
L(e(\lambda_x(u_{(1)}\otimes\omega)))KL(u_{(2)})-
L(e(u_{(1)}\otimes\omega))KL(\lambda_x u_{(2)})=$$
$$=L_xKL(u\otimes\imath(\omega))-KL(u\otimes\imath(\omega))\tau_x+
L_xKL(u_{(1)})q(u_{(2)}\otimes\omega)-
KL(u_{(1)})\tau_x q(u_{(2)}\otimes\omega)+$$
$$+KL(u_{(1)})[\tau_x,q(u_{(2)}\otimes\omega)]-
L(x\otimes e(u_{(1)}\otimes\omega)-e(u_{(1)}\otimes\omega)\otimes x)
KL(u_{(2)})-$$
$$-L(e(u_{(1)}\otimes\omega))L_xKL(u_{(2)})+
L(e(u_{(1)}\otimes\omega))KL(u_{(2)})\tau_x
=L_x Z(u,\omega)-Z(u,\omega)\tau_x.$$
One can write this as
$$Z(x\otimes u,\omega)=L_x Z(u,\omega)-Z(u,\omega)\tau_x-Z(\tau_xu,\omega)-
Z(u,\tau_x\omega).$$
By the induction on the degree of $u,$ we have $Z(u,\omega)=0.$

\section{Geometric interpretation}

Let ${\mathcal M}$ be a smooth manifold,
${\mathfrak F}({\mathcal M})$ be the algebra of
smooth functions on ${\mathcal M}$ and
${\mathcal V}({\mathcal M})=\Der_{\Rs}({\mathfrak F}({\mathcal M}))$
be the Lie algebra of smooth vector fields.
Let us denote ${\mathcal V}={\mathcal V}({\mathcal M})$
and $T({\mathcal V})=T_{\Rs}({\mathcal V}).$

Denote by $\Gamma({\mathcal M},T^n{\mathcal M})$ the space of smooth global
sections of the rank $n$ tensor bundle
$T^n{\mathcal M}=\bigotimes^{n}T{\mathcal M}.$ For example,
$\Gamma({\mathcal M},T^0{\mathcal M})={\mathfrak F}({\mathcal M})$
and
$\Gamma({\mathcal M},T^1{\mathcal M})={\mathcal V}({\mathcal M}).$
Denote
$T({\mathcal M})=\oplus_{n=0}^{\infty}\Gamma({\mathcal M},T^n{\mathcal M}).$
It is well known that
$T({\mathcal M})=T_{{\mathfrak F}({\mathcal M})}({\mathcal V})$
(e.g. [3, Ch. I, Proposition 3.1]). Then there is a natural map
${\mathfrak t}:T({\mathcal V})\to T({\mathcal M}).$

Denote by ${\mathcal D}({\mathcal M})$ the algebra of scalar differential
operators with smooth coefficients on ${\mathcal M}$.
Any vector field is a first order differential operator.
By the definition of $T({\mathcal V}),$
the natural inclusion map $\tau:{\mathcal V}\to {\mathcal D}({\mathcal M})$
can be extended to the algebra homomorphism
$\tau:T({\mathcal V})\to{\mathcal D}({\mathcal M}).$
By the definition of the Lie bracket,
$\tau:x\otimes y-y\otimes x -[x,y]\mapsto 0,\,x,y\in{\mathcal V},$
hence $\tau:I({\mathcal V})\to \{0\}.$

Let us suppose the manifold to be endowed with a smooth affine
connection. Let
$\mu:T({\mathcal M})\to {\mathcal D}({\mathcal M})$ denotes
the ${\mathfrak F}({\mathcal M})$ -- linear
map defined by
$\mu:v_1\otimes\dots\otimes v_n\mapsto\nabla^n_{v_1,\dots,v_n},\,
v_1,\dots,v_n\in{\mathcal V}.$ Here $\nabla^n$ is the n-th order
covariant derivative (in the notation of [3, Ch. III, \S 2]
$\nabla^n_{v_1,\dots,v_n}:f\mapsto f(;v_n;\dots;v_1)$).
For example, if $f,g\in{\mathfrak F}({\mathcal M})$ and
$v\in{\mathcal V}({\mathcal M}),$ then
$\mu(f):g\mapsto fg$ and $\mu(v):g\mapsto v(g).$ The operator
$\mu(v\otimes w)=\nabla^n_{v,w}=vw - (\nabla_v w)$
depends on the connection, as well as the images of the higher degree tensor
fields.

Let $\Sigma({\mathcal M})={\mathfrak t}(\Sigma({\mathcal V}))
\subset T({\mathcal M})$ be the space of (formal sums of)
symmetric tensor fields. By the methods of geometry it may be shown that
there exists a unique map
$\sigma:{\mathcal D}({\mathcal M})\to \Sigma({\mathcal M}),$
such that $\mu\circ\sigma=\id_{{\mathcal D}({\mathcal M})}.$
This map is a surjective ${\mathfrak F}({\mathcal M})$ -- module
homomorphism. The image $\sigma(A)$ is called a symbol of the differential
operator $A.$
In almost the same form the symbol map was introduced
in [6,\S 2] but actually it has been known long before (see the references
in [6]). A proper investigation of the symbol map leads inevitably to the
following natural question: what is the projection
$\Sigma=\sigma\circ\mu:T({\mathcal M})\to\Sigma({\mathcal M})?$
For example, this map is of importance when the symbol of a composition
of two (pseudo)differential operators is considered. In some simple cases,
Sharafutdinov has computed it in the unpublished supplements to [6].

The aforementioned projection $\pi$ is closely related to $\Sigma.$
The space ${\mathcal V}$ may be considered as an algebra with two operations:
the Lie bracket and the covariant derivative (probably Nomizu was the first
who take this view [4, Ch. III, \S 6]). Put
$v\diamond w=\nabla_v w,\,v,w\in{\mathcal V}.$ The corresponding map
$K:T({\mathcal V})\to T({\mathcal V})$ is then connected to $\mu$
by the relation $\mu\circ{\mathfrak t}=\tau\circ K$ [1, Prop. 1], [2, Lemma 2].
This relation is actually a simple consequence of the well known covariant
derivation rules [3, Ch. III, Prop. 2.10].
Let $\pi:T({\mathcal V})\to \Sigma({\mathcal V})$ be the projection
defined above. By the definition,
$\id-\pi:T({\mathcal V})\to K^{-1}(I({\mathcal V})),$
hence $\mu\circ{\mathfrak t}=\mu\circ{\mathfrak t}\circ\pi.$
If ${\mathfrak t}^{-1}:T({\mathcal M})\to T({\mathcal V})$ is any right
inverse of ${\mathfrak t}^{-1},$ then
$\mu=\mu\circ{\mathfrak t}\circ\pi\circ{\mathfrak t}^{-1}.$
The map $\Sigma$ is determined uniquely by the property
$\mu=\mu\circ\Sigma,$
hence $\Sigma={\mathfrak t}\circ\pi\circ{\mathfrak t}^{-1}.$
Note that the equality does not depend on the choice of ${\mathfrak t}^{-1},$
which means
$\pi:\ker({\mathfrak t},T({\mathcal V}))\to
\ker({\mathfrak t},T({\mathcal V})).$
We shall call a linear map with this property
a special one. In other words, the map $T({\mathcal V})\to T({\mathcal V})$
is special if it may be lifted to a ${\mathfrak F}({\mathcal M})$ -
linear map $T({\mathcal M})\to T({\mathcal M}).$

In the expression $\pi=\pi_0(\id+\rho)^{-1}$ the map $\pi_0$ is special,
so it is natural to ask for speciality of $\rho.$
Under natural assumptions on $\imath^{-1}$
this is indeed the case, because the functions $t$ and $r$
are nothing but the
(derivatives of) the torsion tensor and the curvature tensor respectively:
$$t(v\wedge w)=T(v,w),\,t(u\otimes v\wedge w)=(\nabla_u T)(v,w),$$
$$r(v\wedge w)h=R(v,w)h,\,r(u\otimes v\wedge w)h=(\nabla_u R)(v,w)h,$$
etc., where $u,v,w,h\in{\mathcal V}({\mathcal M})$  [4, Ch. III, \S 5], [1,\S 7].
The statement
$R(Q)+\imath(Q)\in K^{-1}(I({\mathcal V})),\,Q\in
\Omega({\mathcal V})\otimes T({\mathcal V})$
may then be considered as a series of commutation relations on the
covariant derivative. For example, if
$Q=u\otimes v\wedge w\otimes h,$ it takes the form
$u\otimes v\otimes w\otimes h-u\otimes w\otimes v\otimes h+
u\otimes t(v\wedge w)\otimes h+t(u\otimes v\wedge w)\otimes h+
u\otimes r(v\wedge w)h+r(u\otimes v\wedge w)h\in K^{-1}(I({\mathcal V})).$
Note that
$\mu\circ{\mathfrak t}=\tau\circ K:K^{-1}(I({\mathcal V}))\to\{0\}.$
Applying this map, we get
$$\nabla^4_{u,v,w,h}-\nabla^4_{u,w,v,h}+\nabla^3_{u,T(v,w),h}+
\nabla^2_{(\nabla_u T)(v,w),h}+\nabla^2_{u,R(v,w)h}+
\nabla_{(\nabla_u R)(v,w)h}=0.$$
Calculating relations of this kind manually is not an easy task
even for relatively small degrees.

\end{document}